\documentclass[12pt]{article}
 \usepackage{latexsym}
 \usepackage{amssymb}
 \usepackage{graphicx}

\newtheorem{Proposition}{Proposition}

\newtheorem{Corollary}{Corollary}

\newtheorem{Algorithm}{Algorithm}

 \newtheorem{theorem}{Theorem}

\newenvironment{example}[1][Example]{\begin{trivlist}
\item[\hskip \labelsep {\bfseries #1}]}{\end{trivlist}}

\newcommand{\J}{{\cal J}}
\newcommand{\A}{{\cal A}}
\newcommand{\B}{{\cal B}}

\newcommand{\I}{{\cal I}}

\newcommand{\qed}{\nobreak \ifvmode \relax \else
      \ifdim\lastskip<1.5em \hskip-\lastskip
      \hskip1.5em plus0em minus0.5em \fi \nobreak
      \vrule height0.75em width0.5em depth0.25em\fi}

\def \ep{\hbox{ }\hfill$\Box$}

\begin{document}
\title{Nonnegative Tensor Factorization, Completely Positive Tensors and an Hierarchical Elimination Algorithm}

\author{
Liqun Qi \thanks{Email: maqilq@polyu.edu.hk. Department of Applied
Mathematics, The Hong Kong Polytechnic University, Hung Hom,
Kowloon, Hong Kong. This author's work was supported by the Hong
Kong Research Grant Council (Grant No. PolyU 501909, 502510, 502111
and 501212).},\ Changqing Xu\thanks{Email: cqxuiit@yahoo.com. School
of Mathematics and Physics, Suzhou University of Science and
Technology, Suzhou, 215009 China.   This author's work was done when
visiting The Hong Kong Polutechnic University and supported by
Jiangsu NSF (No.).}\  and  \ Yi Xu\thanks{Email:
yi.xu1983@gmail.com. Department of Applied Mathematics, The Hong
Kong Polytechnic University, Hung Hom, Kowloon, Hong Kong.}}

\date{\today} \maketitle

\begin{abstract}
\noindent Nonnegative tensor factorization has applications in
statistics, computer vision, exploratory multiway data analysis and
blind source separation. A symmetric nonnegative tensor, which has a
symmetric nonnegative factorization, is called a completely positive
(CP) tensor. The H-eigenvalues of a CP tensor are always
nonnegative. When the order is even, the Z-eigenvalue of a CP tensor
are all nonnegative. When the order is odd, a Z-eigenvector
associated with a positive (negative) Z-eigenvalue of a CP tensor is
always nonnegative (nonpositive).  The entries of a CP tensor obey
some dominance properties.  The CP tensor cone and the copositive
tensor cone of the same order are dual to each other.   We introduce
strongly symmetric tensors and show that a symmetric tensor has a
symmetric binary decomposition if and only if it is strongly
symmetric.  Then we show that a strongly symmetric, hierarchically
dominated nonnegative tensor is a CP tensor, and present a
hierarchical elimination algorithm for checking this.  Numerical
examples are also given. \vspace{3mm}

\noindent {\bf Key words:}\hspace{2mm} completely positive tensor,
eigenvalues, dominance properties, copositive tensor, strongly
symmetric tensor, hierarchical dominance, hierarchical elimination
algorithm. \vspace{3mm}

\noindent {\bf AMS subject classifications (2010):}\hspace{2mm}
15A18; 15A69
  \vspace{3mm}

\end{abstract}


\section{Introduction}
\hspace{4mm}  Nonnegative tensor factorization has applications in
statistics, computer vision, exploratory multiway data analysis and
blind source separation \cite{CZPA, SH}.    As the research topic of
nonnegative matrix factorization is closely related with the theory
of completely positive (CP) matrices \cite{BS, HN, Ka, Xu}, in this
paper, we introduce completely positive (CP) tensors, study their
spectral properties and other properties, demonstrate the dual
relationship between the CP tensor cone and the copositive tensor
cone, and give a checkable sufficient condition for CP tensors by
showing a strongly symmetric, hierarchically dominated nonnegative
tensor is a CP tensor.

Let $\A = \left(a_{i_1 \cdots i_m}\right)$ be a real $m$th order
$n$-dimensional tensor.  Denote the set of all nonnegative vectors
in $\Re^n$ by $\Re^n_+$. For any vector $u \in \Re^n$, $u^m$ is a
rank-one $m$th order symmetric $n$-dimensional tensor $u^m =
\left(u_{i_1}\cdots u_{i_m}\right)$.   If
\begin{equation}\label{e1}
\A = \sum_{k=1}^r \left(u^{(k)}\right)^m,
\end{equation}
where $u^{(k)} \in \Re^n_+$ for $k= 1, \cdots, r$, then $\A$ is called a {\bf completely
positive (CP)} tensor.   The minimum value of $r$ is called the {\bf
CP rank} of $\A$.   The concepts of completely positive tensors and
their CPranks extend the concepts of completely positive matrices
and their CPranks \cite{BS, HN, Ka, Xu}.

The eigenvalues of a CP matrix are always nonnegative.   In the next
section, after summarizing some necessary knowledge about
eigenvalues of tensors, we prove that the H-eigenvalues of a CP
tensor are always nonnegative. We further show that when the order
$m$ is even, the Z-eigenvalue of a CP tensor are all nonnegative,
while when the order $m$ is odd, a Z-eigenvector associated with a
positive (negative) Z-eigenvalue of a CP tensor is always
nonnegative (nonpositive).

In Section 3, we prove some dominance properties which the entries
of a CP tensor must obey.   These properties form some checkable
necessary conditions for a CP tensor.

It is well-known that the CP matrix cone and the copositive matrix
cone are dual to each other.  Recently, motivated by the study of
spectral hypergraph theory, Qi \cite{Qi1} introduced copositive
tensors.   In Section 4, we show that the CP tensor cone and the
copositive tensor cone of the same order are dual to each other.

It is also well-known \cite{Ka} that a diagonally dominated
symmetric nonnegative matrix is a CP matrix.   This forms a
checkable condition for a CP matrix.  To extend this result to CP
tensors, in Section 5, we introduce strongly symmetric tensors and
show that a symmetric tensor is strongly symmetric if and only if it
has a symmetric binary decomposition.   We present a hierarchical
elimination algorithm for checking this.

In Section 6, we further define strongly symmetric, hierarchically
dominated nonnegative tensors and show that a strongly symmetric,
hierarchically dominated nonnegative tensor is a CP tensor.  We show
that the hierarchical elimination algorithm given in Section 5 can
be used to check this condition too.

Some numerical examples are given in Section 7.

Some final remarks are made in Section 8.

For a vector $x \in \Re^n$, denote supp$(x) = \{ i : 1 \le i \le n,
x_i \not = 0 \}$.  For a finite set $S$, $|S|$ denotes its
cardinality.

\section{Eigenvalues of a CP Tensor}

Let $\A = \left(a_{i_1 \cdots i_m}\right)$ be a real $m$th order
$n$-dimensional tensor, and $x \in C^n$. Then
$$\A x^m = \sum_{i_1, \cdots, i_m = 1}^n a_{i_1 \cdots i_m}x_{i_1}
\cdots x_{i_m},$$ and $\A x^{m-1}$ is a vector in $C^n$, with its
$i$th component defined by
$$\left(\A x^{m-1}\right)_i = \sum_{i_2, \cdots, i_m = 1}^n a_{i i_1 \cdots i_m}x_{i_2}
\cdots x_{i_m}.$$ Let $s$ be a positive integer.  Then $x^{[s]}$ is
a vector in $C^n$, with its $i$th component defined by $x_i^s$. We
say that $\A$ is symmetric if its entries $a_{i_1, \cdots, i_m}$ are
invariant for any permutation of the indices.  If $\A x^m \ge 0$ for
all $x \in \Re^n$, then we say that $\A$ is positive semi-definite.
Clearly, only when $m$ is even, a nonzero tensor $A$ can be positive
semi-definite.

The following definitions of eigenvalues, H-eigenvalues,
E-eigenvalues and Z-eigenvalues were introduced in \cite{Qi}.

If $x \in C^n$, $x \not = 0$, $\lambda \in C$,  $x$ and $\lambda$
satisfy
\begin{equation} \label{eig}
\A x^{m-1} = \lambda x^{[m-1]},
\end{equation}
then we call $\lambda$ an {\bf eigenvalue} of $\A$, and $x$ its
corresponding {\bf eigenvector}.   By (\ref{eig}), if $\lambda$ is
an eigenvalue of $\A$ and $x$ is its corresponding eigenvector, then
$$\lambda = {(\A x^{m-1})_j \over x_j^{m-1}},$$
for some $j$ with $x_j \not = 0$.   In particular, if $x$ is real,
then $\lambda$ is also real. In this case, we say that $\lambda$ is
an {\bf H-eigenvalue} of $\A$ and $x$ is its corresponding {\bf
H-eigenvector}.

We say a complex number $\lambda$ is an {\bf E-eigenvalue} of $A$ if
there exists a complex vector $x$ such that
\begin{equation} \label{eig1}
\left\{{Ax^{m-1} = \lambda x, \atop x^Tx = 1.}\right.
\end{equation}
In this case, we say that $x$ is an E-eigenvector of the tensor $A$
associated with the E-eigenvalue $\lambda$.   By (\ref{eig1}), if
$\lambda$ is an E-eigenvalue of $\A$ and $x$ is its E-corresponding
eigenvector, then
$$\lambda = \A x^m.$$
Thus, if $x$ is real, then $\lambda$ is also real. In this case, we
say that $\lambda$ is an {\bf Z-eigenvalue} of $\A$ and $x$ is its
corresponding {\bf Z-eigenvector}.

By \cite{Qi}, we have the following proposition.

\begin{Proposition} \label{p1}   A real $m$th order
$n$-dimensional tensor $A$ has always Z-eigenvalues. If $m$ is even,
then $\A$ always has at least one H-eigenvalue. When $m$ is even,
$\A$ is positive semi-definite if and only if all of its
H-eigenvalues (Z-eigenvalues) are nonnegative.
\end{Proposition}

If the entries of a $m$th order $n$-dimensional tensor $A =
(a_{i_1\cdots i_m})$ are invariant under any permutation of their
indices, then we say that $\A$ is symmetric.  If all the entries of
$\A$ are nonnegative, then we say that $\A$ is a nonnegative tensor.
By (\ref{e1}), a CP tensor is a symmetric nonnegative tensor.  By
\cite{YY}, a symmetric nonnegative tensor has at least one
H-eigenvalue, which is the largest modulus of its eigenvalues.

We now have the following theorem on H-eigenvalues of a CP tensor.

\begin{theorem} \label{t1}
Suppose that $\A = (a_{i_1\cdots i_m})$ is an $m$th order
$n$-dimensional CP tensor, expressed by (\ref{e1}), with $m \ge 2$.
 Then the H-eigenvalues of $\A$ are always nonnegative.
\end{theorem}

\noindent {\bf Proof.} First, assume that $m$ is even. For any $x
\in \Re^n$, we have
$$\A x^m = \sum_{k=1}^r \left(u^{(k)}\right)^mx^m = \sum_{k=1}^r
\left[\left(u^{(k)}\right)^\top x\right]^m \ge 0.$$ Thus, $\A$ is
positive semi-definite.  By Proposition \ref{p1}, all of the
H-eigenvalues are nonnegative.

Now assume that $m$ is odd.   By the discussion before this theorem,
$\A$ has at least one H-eigenvalue.  Suppose that $\lambda$ is an
H-eigenvalue of $\A$, with an H-eigenvector $x$.  Then $x \in \Re^n,
x \not = 0$. By the definition of H-eigenvalue and H-eigenvector, we
have
$$\lambda x^{[m-1]} = \A x^{m-1} = \sum_{k=1}^r
\left(u^{(k)}\right)^mx^{m-1} = \sum_{k=1}^r
\left[\left(u^{(k)}\right)^\top x\right]^{m-1}u^{(k)} \ge 0.$$ Thus,
$\lambda \ge 0$.   This completes the proof.
 \ep

 By (\ref{eig1}), when $m$ is odd, if $\lambda$ is a Z-eigenvalue of a tensor
 $\A$ with a Z-eigenvector $x$, then $-\lambda$ is a Z-eigenvalue of a tensor
 $\A$ with a Z-eigenvector $-x$.  Hence, when $m$ is odd, we cannot
 expect that the Z-eigenvalues of a CP tensor are always
 nonnegative.   However, in this case, we may get strong properties
 of Z-eigenvectors.

\begin{theorem} \label{t2}
Suppose that $\A = (a_{i_1\cdots i_m})$ is an $m$th order
$n$-dimensional CP tensor, expressed by (\ref{e1}), with $m \ge 2$.
When the order $m$ is even, the Z-eigenvalue of a CP tensor are all
nonnegative. When the order $m$ is odd, a Z-eigenvector associated
with a positive (negative) Z-eigenvalue of a CP tensor is always
nonnegative (nonpositive).
\end{theorem}

\noindent {\bf Proof.} The proof of the case that $m$ is even is
similar to the first part of the proof of Theorem \ref{t1}.

Now assume that $m$ is odd.   Suppose that $\lambda$ is a
Z-eigenvalue of $\A$, with an Z-eigenvector $x$.  By the definition
of Z-eigenvalue and Z-eigenvector, we have
$$\lambda x = \A x^{m-1} = \sum_{k=1}^r
\left(u^{(k)}\right)^mx^{m-1} = \sum_{k=1}^r
\left[\left(u^{(k)}\right)^\top x\right]^{m-1}u^{(k)} \ge 0.$$ Thus,
if $\lambda > 0$, then $x \ge 0$, and if $\lambda < 0$, then $x \le
0$. This completes the proof.
 \ep

\section{Dominance Properties of a CP Tensor}

Denote $\I = \{ (i_1,\cdots, i_m) : 1 \le i_k \le n, k = 1, \cdots,
m \}$.  For $(i_1,\cdots, i_m) \in \I$, let $[(i_1,\cdots, i_m)]$ be
the set of all the distinct members in $\{ i_1, \cdots, i_m \}$. For
example, $[(1, 1, 4, 5)] = \{ 1, 4, 5 \}$.

Let $(i_1,\cdots, i_m), (j_1,\cdots, j_m) \in \I$.  We say that
$(i_1,\cdots, i_m)$ is dominated by $(j_1,\cdots, j_m)$, and denote
$(i_1,\cdots, i_m) \preceq (j_1,\cdots, j_m)$ if $[(i_1,\cdots,
i_m)] \subseteq [(j_1,\cdots, j_m)]$.    We say that $(i_1,\cdots,
i_m)$ is similar to $(j_1,\cdots, j_m)$, and denote $(i_1,\cdots,
i_m) \sim (j_1,\cdots, j_m)$ if $[(i_1,\cdots, i_m)] = [(j_1,\cdots,
j_m)]$.

We have the following dominance property for a CP tensor.

\begin{theorem} \label{t3}
Suppose that $\A = (a_{i_1\cdots i_m})$ is an $m$th order
$n$-dimensional CP tensor, expressed by (\ref{e1}), with $m \ge 2$.
If $(i_1,\cdots, i_m) \preceq (j_1,\cdots, j_m)$ and $a_{j_1\cdots j_m} \not = 0$, then $a_{i_1\cdots i_m} > 0$.
\end{theorem}

\noindent {\bf Proof.} We have that
$$0 \not = a_{j_1,\cdots, j_m} = \sum_{k=1}^r u^{(k)}_{j_1}\cdots u^{(k)}_{j_m}.$$
Since $u^{(k)} \in \Re^n_+$ for $k= 1, \cdots, r$, at least for one $k = \bar k$, $u^{(\bar k)}_{j_1} > 0, \cdots, u^{(\bar k)}_{j_m} > 0.$  Since $\{ i_1,\cdots, i_m \} \subseteq \{ j_1,\cdots, j_m \}$, this implies that
 $u^{(\bar k)}_{i_1} > 0, \cdots, u^{(\bar k)}_{i_m} > 0.$   Therefore,
 $$a_{i_1,\cdots, i_m} = \sum_{k=1}^r u^{(k)}_{i_1}\cdots u^{(k)}_{i_m} > 0.$$
 This completes the proof.  \ep

\begin{Corollary} \label{c1}
Suppose that $\A = (a_{i_1\cdots i_m})$ is an $m$th order
$n$-dimensional CP tensor, expressed by (\ref{e1}), with $m \ge 2$.
If $(i_1,\cdots, i_m) \sim (j_1,\cdots, j_m)$, then $a_{j_1\cdots
j_m} = 0$ if and only if $a_{i_1\cdots i_m} = 0$.
\end{Corollary}

When $m=2$, this property can be derived from the symmetric property
of the matrix $\A$.   When $m > 2$, this property cannot be derived
from the symmetric property of the tensor $\A$.   For example, for a
third order CP tensor $\A = (a_{ijk})$, we have $a_{iij} = a_{ijj}$
for all $i$ and $j$, satisfying $1 \le i, j \le n$.  But this is not
true for a general third order symmetric tensor.   This motivates us
to introduce strongly symmetric tensors in Section 5.

 Suppose that $(j_1, \cdots, j_m) \in \I$ and $I = \left\{ \left(i^{(1)}_1,\cdots, i^{(1)}_m\right), \cdots, \left(i^{(s)}_1,\cdots, i^{(s)}_m\right) \right\} \subseteq \I$. Assume that
 $(i^{(p)}_1,\cdots, i^{(p)}_m) \preceq (j_1, \cdots, j_m)$ for $p = 1, \cdots, s$, and for any index $i \in \{ j_1, \cdots, j_m \}$, if it appears $t$ times in $\{ j_1, \cdots, j_m \}$, then it appears in $I$ $st$ times.   Then we call $I$ an {\bf $s$-duplicate} of $(j_1, \cdots, j_m)$.

 We have the following strong dominance property for a CP tensor.

\begin{theorem} \label{t4}
Suppose that $\A = (a_{i_1\cdots i_m})$ is an $m$th order
$n$-dimensional CP tensor, expressed by (\ref{e1}), with $m \ge 2$.
Assume that $I = \left\{ \left(i^{(1)}_1,\cdots, i^{(1)}_m\right),
\cdots, \left(i^{(s)}_1,\cdots, i^{(s)}_m\right) \right\}$ is an
s-duplicate of $(j_1,\cdots j_m) \in \I$.  Then
$${1 \over s}\sum_{p=1}^s a_{i^{(p)}_1\cdots i^{(p)}_m} \ge a_{j_1\cdots j_m}.$$
\end{theorem}

\noindent {\bf Proof.} We have that $${1 \over s}\sum_{p=1}^s
a_{i^{(p)}_1\cdots i^{(p)}_m} = \sum_{k=1}^r {1 \over s}\sum_{p=1}^s
u^{(k)}_{i^{(p)}_1}\cdots u^{(k)}_{i^{(p)}_m} \ge \sum_{k=1}^r
\left(\prod_{p=1}^s u^{(k)}_{i^{(p)}_1}\cdots u^{(k)}_{i^{(p)}_m}
\right)^{1 \over s} = \sum_{k=1}^r u^{(k)}_{j_1}\cdots u^{(k)}_{j_m}
= a_{j_1\cdots j_m},$$ where the inequality is due to the fact that
the geometric mean of some positive numbers is never greater than
their arithmetic mean.    This completes the proof. \ep

\begin{Corollary} \label{c2}
Suppose that $\A = (a_{i_1\cdots i_m})$ is an $m$th order
$n$-dimensional CP tensor, expressed by (\ref{e1}), with $m \ge 2$.
Assume that $(j_1,\cdots j_m) \in \I$.  Then
$${1 \over m}\sum_{p=1}^m a_{j_p\cdots j_p} \ge a_{j_1\cdots j_m}.$$
\end{Corollary}

\section{The CP Tensor Cone and the Copositive Tensor Cone}

Denote the set of all $m$th order $n$-dimensional CP tensors by
$CP_{m,n}$.   By (\ref{e1}), it is easy to see that $CP_{m,n}$ is a
closed convex cone.   Suppose that $\B$ is a real $m$th order $n$th
dimensional symmetric tensor.  If for all $x \in \Re^n_+$, we have
$\B x^m \ge 0$, then $\B$ is called a {\bf copositive tensor}
\cite{Qi1}.  Denote the set of all $m$th order $n$-dimensional
copositive tensors by $COP_{m,n}$.  Then, it is also easy to see
that $COP_{m,n}$ is a closed convex cone.    When $m =2$, a
classical result is that the CP matrix cone and the copositive
matrix cone are dual to each other.  We now extend this result to
the CP tensor cone and the copositive tensor cone.

Let $\A = (a_{i_1\cdots i_m})$ and $\B = (b_{i_1\cdots i_m})$ be two
real $m$th order $n$-dimensional symmetric tensors.  Their inner
product is defined as
$$\A \bullet \B = \sum_{i_1,\cdots, i_m=1}^n a_{i_1\cdots i_m}b_{i_1\cdots
i_m}.$$

\begin{theorem} \label{t0}
Let $m \ge 2$ and $n \ge 1$.   Then $CP_{m,n}$ and $COP_{m,n}$ are
dual to each other.
\end{theorem}

\noindent {\bf Proof.}  Suppose that $\B$ is an $m$th order
$n$-dimensional copositive tensor.  For any $\A \in CP_{m, n}$, by
definition, we may assume that $\A$ can be expressed by (\ref{e1}).
Since $\B$ is a copositive tensor, by definition, $\B
\left(u^{(k)}\right)^m \ge 0$, for $k = 1, \cdots, r$.   Thus,
$$\A \bullet \B = \sum_{k=1}^r \B
\left(u^{(k)}\right)^m \ge 0.$$ Thus, $\B$ is in the dual cone of
$CP_{m,n}$.

On the other hand, assume that $\B$ is in the dual cone of
$CP_{m,n}$.   Let $x \in \Re^n_+$.  Then $x^m$ is an $m$ order
$n$-dimensional CP tensor, i.e., $x \in CP_{m, n}$.  We have $\B x^m
= \B \bullet x^m \ge 0$.  This shows that $\B$ is a copositive
tensor.

Together, we see that $CP_{m,n}$ and $COP_{m,n}$ are dual to each
other.
 \ep

\section{Strongly Symmetric Tensors}

Suppose that $\A = \left(a_{i_1 \cdots i_m}\right)$ is a real $m$th
order $n$-dimensional tensor.  If for any $(i_1,\cdots, i_m) \sim
(j_1,\cdots, j_m)$, $(i_1,\cdots, i_m), (j_1,\cdots, j_m) \in \I$,
we have $a_{i_1\cdots i_m} = a_{j_1\cdots j_m}$, then we say that
$\A$ is a {\bf strongly symmetric} tensor.   Clearly, a strongly
symmetric tensor is a symmetric tensor.  It is also clear that a
linear combination of strongly symmetric tensors is still a strongly
symmetric tensor.   Thus, the set of all real $m$th order
$n$-dimensional strongly symmetric tensors is a linear space.

Let $\A = \left(a_{i_1 \cdots i_m}\right)$ be a real $m$th order
$n$-dimensional symmetric tensor.   If
\begin{equation}\label{e4}
\A = \sum_{k=1}^r \alpha_k\left(v^{(k)}\right)^m,
\end{equation}
where $\alpha_k$ are real numbers and $v^{(k)}$ are binary vectors
in $\Re^n$ for $k= 1, \cdots, r$, then we say that $\A$ has a {\bf
symmetric binary decomposition}, which is not a nonnegative tensor
factorization, but a general symmetric tensor decomposition
\cite{CGLM, KB}.

It is easy to show the following proposition.

\begin{Proposition} \label{p2}
Suppose that $\A = (a_{i_1\cdots i_m})$ is a real $m$th order
$n$-dimensional tensor with a symmetric binary decomposition. Then
$\A$ is strongly symmetric.
\end{Proposition}

\noindent {\bf Proof.}  Suppose that $\A = (a_{i_1\cdots i_m})$ is
expressed by (\ref{e4}).   Assume that $(i_1,\cdots, i_m) \sim
(j_1,\cdots, j_m)$.  Then
$$a_{i_1\cdots i_m} = \sum \left\{ \alpha_k : (i_1, \cdots i_m) \preceq {\rm
supp}\left(v^{(k)}\right) \right\} = \sum \left\{ \alpha_k : (j_1,
\cdots j_m) \preceq {\rm supp}\left(v^{(k)}\right) \right\} =
a_{j_1\cdots j_m}.$$ This completes the proof.  \ep

\medskip

For $k = 1, \cdots, m$, let
$$\I_k = \left\{ (i_1, \cdots, i_m) \in \I : \left|[(i_1, \cdots, i_m)]\right| = k
\right\}.$$ Then $\I_1, \cdots, \I_m$ form a partition of $\I$.

For $k = 1, \cdots, m$, let
$$\I_{k+} = \left\{ (i_1, \cdots, i_k, i_k, \cdots, i_k) \in \I_k : 1 \le i_1
< i_2 \cdots < i_k \le n \right\}.$$ Then $\I_{k+}$ is the
``representative'' set of $\I_k$ in the sense that any member in
$\I_k$ is similar to a member of $\I_{k+}$ and no two members in
$\I_{k+}$ are similar.

Suppose that $\A = \left(a_{i_1 \cdots i_m}\right)$ is a real $m$th
order $n$-dimensional tensor.  For $k = 1, \cdots, m$, let
$$\I_{k+}(\A) = \left\{ (i_1, \cdots, i_k, i_k, \cdots, i_k) \in \I_{k+} : a_{i_1\cdots i_ki_k\cdots i_k} \not = 0 \right\}.$$

\medskip

We now construct a hierarchical elimination algorithm to obtain
symmetric binary decomposition of a strongly symmetric tensor.
Suppose that $\A = \left(a_{i_1 \cdots i_m}\right)$ is a real $m$th
order $n$-dimensional strongly symmetric tensor.

\begin{Algorithm} \label{a1}

Step 0. Let $k = 0$ and $\A^{(0)}=\left(a^{(0)}_{i_1\cdots
i_m}\right)$ be defined by $\A^{(0)} = \A$.

Step 1. For any $e = (i_1, \cdots, i_{m-k}, \cdots,
 i_{m-k}) \in \I_{(m-k)+}\left(\A^{(k)}\right)$, let $v^e \in \Re^n_+$
 be a binary vector such that $v^e_{i_1} = \cdots  = v^e_{i_{m-k}}
 = 1$ and $v^e_i = 0$ if $i \not \in \{ i_1, \cdots, i_{m-k} \}$.

Let $\A^{(k+1)}=\left(a^{(k+1)}_{i_1\cdots i_m}\right)$ be defined
by
\begin{equation} \label{e5}
 \A^{(k+1)} = \A^{(k)} - \sum \left\{ a^{(k)}_{i_1\cdots i_{m-k}\cdots
 i_{m-k}} \left(v^{e}\right)^m : e = (i_1, \cdots, i_{m-k}, \cdots,
 i_{m-k}) \in \I_{(m-k)+}\left(\A^{(k)}\right) \right\}.
\end{equation}

Step 2. Let $k = k+1$.  If $k = m$, stop.  Otherwise, go to Step 1.
\end{Algorithm}

\begin{theorem} \label{t5}
Suppose that $\A = (a_{i_1\cdots i_m})$ is a real $m$th order
$n$-dimensional strongly symmetric tensor.   Then we have $\A^{(m)}
= 0$ in Algorithm \ref{a1}, i.e., we have
\begin{equation} \label{e6}
 \A = \sum_{k=0}^{m-1} \sum \left\{ a^{(k)}_{i_1\cdots i_{m-k}\cdots
 i_{m-k}} \left(v^{e}\right)^m : e = (i_1, \cdots, i_{m-k}, \cdots,
 i_{m-k}) \in \I_{(m-k)+}\left(\A^{(k)}\right) \right\}.
\end{equation}

Thus, a symmetric tensor has a symmetric binary decomposition if and
only if it is strongly symmetric.
\end{theorem}

\noindent {\bf Proof.} For $k=1,\cdots,m$, we now show by induction
that $\A^{(k)}$ is strongly symmetric, and $\I_{(m-p)+}(\A^{(k)}) =
\emptyset$ for $p = 0,\cdots, k-1$.

By Step 0 and the assumption, $\A^{(0)}$ is strongly symmetric.

For $k=0,\cdots, m-1$, assume that $\A^{(k)}$ is strongly symmetric,
and $\I_{(m-p)+}(\A^{(k)}) = \emptyset$ for $p = 0,\cdots, k-1$ if
$k \ge 1$.  By (\ref{e5}) and Proposition \ref{p2}, $\A^{(k+1)}$ is
a linear combination of strongly symmetric tensors, thus also a
strongly symmetric tensor.   As in this iteration $|$supp$(v^e)|=
m-k$ for all $v^e$ in (\ref{e5}), $a^{(k+1)}_{i_1\cdots i_m} =
a^{(k)}_{i_1\cdots i_m} = 0$ if $|[(i_1,\cdots, i_m)]| > m-k$. Thus,
$\I_{(m-p)+}(\A^{(k+1)}) = \emptyset$ for $p = 0,\cdots, k-1$.  By
(\ref{e5}), we also have $\I_{(m-k)+}(\A^{(k+1)}) = \emptyset$. The
induction proof is completed.

This shows that $\A^{(m)} = 0$. By this and (\ref{e5}), we have
(\ref{e6}).   Thus, a strongly symmetric tensor has a symmetric
binary decomposition.    By this and Proposition \ref{p2}, the last
conclusion also holds. This completes the proof. \ep

\section{Strongly Symmetric, Hierarchically Dominated Tensors}

In (\ref{e6}), if all the coefficients $a^{(k)}_{i_1\cdots
i_{m-k}\cdots i_{m-k}}$ are nonnegative, then $\A$ is a CP tensor.
In this section, we explore a sufficient condition for this.

For $p=1,\cdots, m-1$, and $q=1, \cdots, m-p$, for any $(i_1,\cdots,
i_p,i_p,\cdots, i_p) \in \I_{p+}$, define
$$\J_q(i_1,\cdots, i_p) = \left\{ (j_1,\cdots, j_{p+q},\cdots, j_{p+q})
\in \I_{(p+q)+} : (i_1,\cdots, i_p,\cdots, i_p) \preceq (j_1,\cdots,
j_{p+q},\cdots, j_{p+q}) \right\}.$$

An $m$th order $n$-dimensional strongly symmetric nonnegative tensor
$\A = (a_{i_1\cdots i_m})$ is said to be {\bf hierarchically
dominated} if for $p=1,\cdots, m-1$, and any $(i_1,\cdots,
i_p,i_p,\cdots, i_p) \in \I_{p+}$, we have
\begin{equation} \label{e7}
a_{i_1\cdots i_pi_p\cdots i_p} \ge \sum \left\{ a_{j_1\cdots
j_{p+1}\cdots j_{p+1}} : (j_1,\cdots, j_{p+1},j_{p+1},\cdots,
j_{p+1}) \in \J_1(i_1,\cdots, i_p)\right\}.
\end{equation}

Suppose that $\A$ is an $m$th order $n$-dimensional strongly
symmetric, hierarchically dominated nonnegative tensor.  By
(\ref{e7}), for $p=1,\cdots, m-2$, and any $(i_1,\cdots,
i_p,i_p,\cdots, i_p) \in \I_{p+}$, we have $$\begin{array} {rcl}
a_{i_1\cdots i_pi_p\cdots i_p} & \ge &\sum \left\{ a_{j_1\cdots
j_{p+1}\cdots j_{p+1}} : (j_1,\cdots, j_{p+1},j_{p+1},\cdots,
j_{p+1}) \in \J_1(i_1,\cdots, i_p)\right\}\\
& \ge & \sum \left\{ \sum \left\{ a_{l_1\cdots l_{p+2}\cdots
j_{p+2}} : (l_1,\cdots, l_{p+2},\cdots, l_{p+2}) \in
\J_1(j_1,\cdots, j_{p+1})\right\} \right.\\
&& \left.: (j_1,\cdots, j_{p+1},\cdots, j_{p+1}) \in
\J_1(i_1,\cdots, i_p)\right\}\\
&\ge & \sum \left\{ a_{l_1\cdots l_{p+2}\cdots j_{p+2}} :
(l_1,\cdots, l_{p+2},\cdots, l_{p+2}) \in \J_2(i_1,\cdots,
i_p)\right\}.
\end{array}$$
Thus, by induction, we may prove the following proposition.

\begin{Proposition} \label{p3}
Suppose that $\A = (a_{i_1\cdots i_m})$ is an $m$th order
$n$-dimensional strongly symmetric, hierarchically dominated
nonnegative tensor.  Then for $p=1,\cdots, m-1$, and $q=1, \cdots,
m-p$, for any $(i_1,\cdots, i_p,i_p,\cdots, i_p) \in \I_{p+}$, we
have
\begin{equation} \label{e8}
a_{i_1\cdots i_pi_p\cdots i_p} \ge \sum \left\{ a_{j_1\cdots
j_{p+q}\cdots j_{p+q}} : (j_1,\cdots, j_{p+q},\cdots, j_{p+q}) \in
\J_q(i_1,\cdots, i_p)\right\}.
\end{equation}
\end{Proposition}

With this proposition, we may prove the following main theorem of
this section.

\begin{theorem} \label{t6}
Suppose that $\A = (a_{i_1\cdots i_m})$ is an $m$th order
$n$-dimensional strongly symmetric, hierarchically dominated
nonnegative tensor. Then $\A^{(k)}$ are nonnegative for $k=0,\cdots,
m-1$, in Algorithm \ref{a1}.   Thus, $\A$ is a CP tensor.  Thus, a
strongly symmetric, hierarchically dominated nonnegative tensor is a
CP tensor.
\end{theorem}

\noindent {\bf Proof.} For $k=1,\cdots,m-1$, we now show by
induction that $\A^{(k)}$ is a strongly symmetric, hierarchically
dominated nonnegative tensor.

By Step 0 and the assumption, $\A^{(0)}$ is a strongly symmetric,
hierarchically dominated nonnegative tensor.

For $k=0,\cdots, m-1$, assume that $\A^{(k)}$ is a strongly
symmetric, hierarchically dominated nonnegative tensor.  We now
consider $\A^{(k+1)}$.

By the proof of Theorem \ref{t5}, $\A^{(k+1)}$ is also strongly
symmetric and for $p=0,\cdots, k$, and any $(i_1,\cdots, i_m) \in
\I_{(m-p)+}$, $a^{(k+1)}_{i_1\cdots i_m} = 0$.  By strong symmetry
of $\A^{(k+1)}$, for $p=0,\cdots, k$, and any $(i_1,\cdots, i_m) \in
\I_{m-p}$, $a^{(k+1)}_{i_1\cdots i_m} = 0$.

Now for $p= k+1,\cdots, m-1$, and any $(i_1,\cdots, i_m) \in
\I_{(m-p)+}$, by (\ref{e5}),
\begin{eqnarray} \label{e9}
 && a^{(k+1)}_{i_1\cdots i_{m-p}\cdots i_{m-p}} \nonumber \\
& =  &a^{(k)}_{i_1\cdots i_{m-p}\cdots i_{m-p}}  \nonumber \\ &- &
\sum \left\{ a^{(k)}_{l_1\cdots l_{m-k}\cdots
 l_{m-k}} : (l_1, \cdots, l_{m-k}, \cdots,
 l_{m-k}) \in \J_{p-k}(i_1,\cdots, i_{m-p}) \right\}.
\end{eqnarray}
By Proposition \ref{p3}, the right hand side of (\ref{e6}) is
nonnegative.   Thus, for $p= k+1,\cdots, m-1$, and any $(i_1,\cdots,
i_m) \in \I_{(m-p)+}$, $a^{(k+1)}_{i_1\cdots i_{m-p}\cdots i_{m-p}}
\ge 0$.  By strong symmetry of $\A^{(k+1)}$, for $p= k+1,\cdots,
m-1$, and any $(i_1,\cdots, i_m) \in \I_{m-p}$,
$a^{(k+1)}_{i_1\cdots i_{m-p}\cdots i_{m-p}} \ge 0$.  This shows
that $\A^{(k+1)}$ is nonnegative.

Since $\A^{(k)} = (a^{(k)}_{i_1\cdots i_m})$ is hierarchically
dominated,  for $p=1,\cdots, m-1$, and any \\ $(i_1,\cdots,
i_p,\cdots, i_p) \in \I_{p+}$, we have
\begin{equation} \label{e10}
a^{(k)}_{i_1\cdots i_p\cdots i_p} \ge \sum \left\{
a^{(k)}_{j_1\cdots j_{p+1}\cdots j_{p+1}} : (j_1,\cdots,
j_{p+1},\cdots, j_{p+1}) \in \J_1(i_1,\cdots, i_p)\right\}.
\end{equation}

By (\ref{e6}), we have
\begin{eqnarray} \label{e11}
&& a^{(k+1)}_{i_1\cdots i_p\cdots i_p} \nonumber \\ & = &
a^{(k)}_{i_1\cdots i_p\cdots i_p} - \sum \left\{ a^{(k)}_{l_1\cdots
l_{m-k}\cdots
 l_{m-k}} : (l_1, \cdots, l_{m-k}, \cdots,
 l_{m-k}) \in \J_{m-p-k}(i_1,\cdots, i_p) \right\}
\end{eqnarray}
and
\begin{eqnarray} \label{e12}
&& a^{(k+1)}_{j_1\cdots j_{p+1}\cdots j_{p+1}} \nonumber \\  &
 =  & a^{(k)}_{j_1\cdots j_{p+1}\cdots
j_{p+1}} \nonumber \\ &-& \sum \left\{ a^{(k)}_{l_1\cdots
l_{m-k}\cdots
 l_{m-k}} : (l_1, \cdots, l_{m-k}, \cdots,
 l_{m-k}) \in \J_{m-p-1-k}(j_1,\cdots, j_{p+1}) \right\}
\end{eqnarray}
Comparing (\ref{e10}), (\ref{e11}) and (\ref{e12}), for $p=1,\cdots,
m-1$, and any \\ $(i_1,\cdots, i_p,\cdots, i_p) \in \I_{p+}$, we
have $$a^{(k+1)}_{i_1\cdots i_p\cdots i_p} \ge \sum \left\{
a^{(k+1)}_{j_1\cdots j_{p+1}\cdots j_{p+1}} : (j_1,\cdots,
j_{p+1},\cdots, j_{p+1}) \in \J_1(i_1,\cdots, i_p)\right\}.$$ Thus,
$\A^{(k+1)}$ is also hierarchically dominated. The induction proof
is completed.

Hence, $\A$ is a CP tensor. Therefore, a strongly symmetric,
hierarchically dominated nonnegative tensor is a CP tensor. This
completes the proof. \ep

When $m=2$, Theorem \ref{t6} implies Kaykobad's result \cite{Ka}.

\begin{Corollary} \label{c3}
Suppose that $\A = (a_{i_1\cdots i_m})$ is a real $m$th order
$n$-dimensional strongly symmetric, hierarchically dominated
nonnegative tensor. Then the CPrank of $\A$ is not bigger than
$\sum_{k=0}^{m-1} \left({n \atop m-k}\right)$.
\end{Corollary}
\noindent {\bf Proof.} By (\ref{e6}), the CPrank of $\A$ is not
bigger than $$\sum_{k=0}^{m-1}
\left|\I_{(m-k)+}\left(\A^{(k)}\right)\right| =  \sum_{k=0}^{m-1}
\left({n \atop m-k}\right).$$  This completes the proof.  \ep

\section{Numerical Examples}
In this section, we present some strongly symmetric, hierarchically
dominated nonnegative tensors with $m=3, n=10$ and $m=4, n=10$, and
use Algorithm \ref{a1} to decompose them.
\begin{example}
$\A$ is a strongly symmetric, hierarchically dominated nonnegative
tensor. The entries of $\A$, whose index sets are not similar to the
index sets of the entries defined below, are zero.

The $m=3$, $n=10$ case:

(1) $\A(1,1,1)=1$, $\A(2,2,2)=5$, $\A(3,3,3)=3$, $\A(4,4,4)=2$, $\A(5,5,5)=4$, $\A(6,6,6)=2$, $\A(7,7,7)=2$, $\A(8,8,8)=2$, $\A(9,9,9)=5$, $\A(10,10,10)=4$, $\A(1,5,5)=1$, $\A(2,3,3)=1$, $\A(2,6,6)=1$, $\A(2,8,8)=1$, $\A(3,4,4)=1$, $\A(3,5,5)=1$,
$\A(4,5,5)=1$, $\A(5,9,9)=1$, $\A(6,9,9)=1$, $\A(7,9,9)=1$, $\A(7,10,10)=1$, $\A(8,10,10)=1$, $\A(9,10,10)=1$, $\A(2,6,9)=1$,
$\A(2,8,10)=1$, $\A(3,4,5)=1$, $\A(7,9,10)=1$.
\\

(2) $\A(1,1,1)=2$, $\A(2,2,2)=5$, $\A(3,3,3)=6$, $\A(4,4,4)=2$, $\A(5,5,5)=3$, $\A(8,8,8)=6$, $\A(9,9,9)=6$, $\A(10,10,10)=4$, $\A(1,5,5)=1$, $\A(1,10,10)=1$, $\A(2,3,3)=1$, $\A(2,8,8)=1$, $\A(2,9,9)=2$, $\A(2,10,10)=1$, $\A(3,4,4)=1$, $\A(3,8,8)=2$, $\A(3,9,9)=2$, $\A(4,8,8)=1$, $\A(5,8,8)=1$, $\A(5,10,10)=1$, $\A(8,9,9)=1$, $\A(9,10,10)=1$, $\A(1,5,10)=1$, $\A(2,3,9)=1$, $\A(2,9,10)=1$, $\A(3,4,8)=1$, $\A(3,8,9)=1$.
\\

(3)  $\A(2,2,2)=4$, $\A(3,3,3)=6$, $\A(4,4,4)=7$, $\A(5,5,5)=4$, $\A(7,7,7)=4$, $\A(8,8,8)=6$, $\A(9,9,9)=4$, $\A(10,10,10)=3$, $\A(2,3,3)=1$, $\A(2,4,4)=1$, $\A(2,5,5)=1$, $\A(2,8,8)=1$, $\A(3,4,4)=1$, $\A(3,5,5)=1$, $\A(3,7,7)=1$, $\A(3,8,8)=2$, $\A(4,5,5)=2$, $\A(4,7,7)=1$, $\A(4,9,9)=1$, $\A(4,10,10)=1$, $\A(7,8,8)=1$, $\A(7,9,9)=1$, $\A(8,9,9)=1$, $\A(8,10,10)=1$, $\A(9,10,10)=1$, $\A(2,3,8)=1$, $\A(2,4,5)=1$, $\A(3,4,5)=1$, $\A(3,7,8)=1$, $\A(4,7,9)=1$, $\A(8,9,10)=1$.
\\

The $m=4, n=10$ case:

(1) $\A(1,1,1,1)=1$, $\A(2,2,2,2)=6$, $\A(4,4,4,4)=6$, $\A(5,5,5,5)=2$, $\A(6,6,6,6)=3$, $\A(7,7,7,7)=4$, $\A(8,8,8,8)=8$, $\A(9,9,9,9)=12$, $\A(10,10,10,10)=4$, $\A(1,10,10,10)=1$, $\A(2,4,4,4)=2$, $\A(2,8,8,8)=2$, $\A(2,9,9,9)=2$, $\A(4,8,8,8)=2$, $\A(4,9,9,9)=2$, $\A(5,7,7,7)=1$, $\A(5,9,9,9)=1$, $\A(6,7,7,7)=1$, $\A(6,9,9,9)=1$, $\A(6,10,10,10)=1$, $\A(7,9,9,9)=2$, $\A(8,9,9,9)=3$, $\A(8,10,10,10)=1$, $\A(9,10,10,10)=1$, $\A(2,4,8,8)=1$, $\A(2,4,9,9)=1$, $\A(2,8,9,9)=1$, $\A(4,8,9,9)=1$, $\A(5,7,9,9)=1$, $\A(6,7,9,9)=1$, $\A(8,9,10,10)=1$, $\A(2,4,8,9)=1$.
\\

(2)  $\A(1,1,1,1)=9$, $\A(2,2,2,2)=6$, $\A(3,3,3,3)=8$, $\A(4,4,4,4)=1$, $\A(5,5,5,5)=1$, $\A(6,6,6,6)=4$, $\A(7,7,7,7)=6$, $\A(8,8,8,8)=6$, $\A(9,9,9,9)=9$, $\A(10,10,10,10)=2$, $\A(1,2,2,2)=1$, $\A(1,3,3,3)=2$, $\A(1,5,5,5)=1$, $\A(1,7,7,7)=1$, $\A(1,8,8,8)=2$, $\A(1,9,9,9)=2$, $\A(2,3,3,3)=1$, $\A(2,6,6,6)=2$, $\A(2,7,7,7)=2$, $\A(3,6,6,6)=1$, $\A(3,8,8,8)=2$, $\A(3,9,9,9)=2$, $\A(4,9,9,9)=1$, $\A(6,7,7,7)=1$, $\A(7,9,9,9)=1$, $\A(7,10,10,10)=1$, $\A(8,9,9,9)=2$, $\A(9,10,10,10)=1$, $\A(1,2,7,7)=1$, $\A(1,3,8,8)=1$, $\A(1,3,9,9)=1$, $\A(1,8,9,9)=1$, $\A(2,3,6,6)=1$, $\A(2,6,7,7)=1$, $\A(3,8,9,9)=1$, $\A(7,9,10,10)=1$, $\A(1,3,8,9)=1$.
\\

(3) $\A(1,1,1,1)=18$, $\A(2,2,2,2)=6$, $\A(3,3,3,3)=4$, $\A(4,4,4,4)=2$, $\A(5,5,5,5)=26$, $\A(6,6,6,6)=18$, $\A(8,8,8,8)=8$, $\A(9,9,9,9)=24$, $\A(10,10,10,10)=8$, $\A(1,5,5,5)=6$, $\A(1,6,6,6)=4$, $\A(1,8,8,8)=2$, $\A(1,9,9,9)=4$, $\A(1,10,10,10)=2$, $\A(2,5,5,5)=2$, $\A(2,6,6,6)=2$, $\A(2,9,9,9)=2$, $\A(3,4,4,4)=1$, $\A(3,9,9,9)=2$, $\A(3,10,10,10)=1$, $\A(4,9,9,9)=1$, $\A(5,6,6,6)=6$, $\A(5,8,8,8)=3$, $\A(5,9,9,9)=7$, $\A(5,10,10,10)=2$, $\A(6,8,8,8)=2$, $\A(6,9,9,9)=4$, $\A(8,9,9,9)=1$, $\A(9,10,10,10)=3$, $\A(1,5,6,6)=2$, $\A(1,5,8,8)=1$, $\A(1,5,9,9)=2$, $\A(1,5,10,10)=1$, $\A(1,6,8,8)=1$, $\A(1,6,9,9)=1$, $\A(1,9,10,10)=1$, $\A(2,5,6,6)=1$, $\A(2,5,9,9)=1$, $\A(2,6,9,9)=1$, $\A(3,4,9,9)=1$, $\A(3,9,10,10)=1$, $\A(5,6,8,8)=1$, $\A(5,6,9,9)=2$, $\A(5,8,9,9)=1$, $\A(5,9,10,10)=1$, $\A(1,5,6,8)=1$, $\A(1,5,6,9)=1$, $\A(1,5,9,10)=1$, $\A(2,5,6,9)=1$.

The vectors decomposed from Algorithm \ref{a1} are shown in Tables
\ref{t1},\ref{t2},\ref{t3},\ref{t4},\ref{t5} and \ref{t6}, in which
$v$ rows are the nonzero values of vectors, $p$ rows are the indices
of the nonzero values.
\begin{table}
\begin{center}
\begin{tabular}{|c|c|c|c|c|c|c|c|c|c|c|c|c|c|c|c|c|c|c|c|c|c|c|c|c|c|c|c|c|c|c|c|c|}
  \hline
  v&1&1&1&1&1&1&1&1.2599&1&1&1 \\
  \hline
    p& 2  6  9&2   8  10&3  4  5&7   9  10&1  5&2  3&5  9&2&3&4&5 \\
  \hline
  v&1&1&1&1.2599&1.2599 & & & & & &\\
  \hline
  p&6&7&8&9&10 & & & & & &\\
  \hline
\end{tabular}
\caption{$n=10,m=3$ (1)}\label{t1}
\end{center}
\end{table}

\begin{table}
\begin{center}
\begin{tabular}{|c|c|c|c|c|c|c|c|c|c|c|c|c|c|c|c|c|c|c|c|c|c|c|c|c|c|c|c|c|c|c|c|c|}
  \hline
  v&1&1&1&1&1&1&1&1&1.2599&1.4422 \\
  \hline
    p&1   5  10&2  3  9&2   9  10&3  4  8&3  8  9&2  8&5  8&1&2&3 \\
  \hline
  v&1&1&1.2599&1.4422&1.2599& & & & & \\
  \hline
  p&4&5&8&9&10& & & & & \\
  \hline
\end{tabular}
\caption{$n=10,m=3$ (2)}\label{t2}
\end{center}

\end{table}

\begin{table}
\begin{center}
\begin{tabular}{|c|c|c|c|c|c|c|c|c|c|c|c|c|c|c|c|c|c|c|c|c|c|c|c|c|c|c|c|c|c|c|c|c|}
  \hline
  v&1&1&1&1&1&1&1&1.2599&1.4422 \\
  \hline
    p&2  3  8&2  4  5&3  4  5&3  7  8&4  7  9&8   9  10&4  10&2&3 \\
  \hline
  v&1.4422&1.2599&1.2599&1.4422&1.2599&1& & &\\
  \hline
  p&4&5&7&8&9&10& & &\\
  \hline
\end{tabular}
\caption{$n=10,m=3$ (3)}\label{t3}
\end{center}

\end{table}

\begin{table}
\begin{center}
\begin{tabular}{|c|c|c|c|c|c|c|c|c|c|c|c|c|c|c|c|c|c|c|c|c|c|c|c|c|c|c|c|c|c|c|c|c|}
  \hline
  v&1&1&1&1&1&1&1&1&1&1&1 \\
  \hline
    p&2  4  8  9&5  7  9&6  7  9&8   9  10&1  10&2  4&2  8&2  9&4  8&4  9&6  10\\
  \hline
  v&1&1.1892&1.1892&1&1&1.1892&1.3161&1.4953&1& &\\
  \hline
  p&8  9&2&4&5&6&7&8&9&10& & \\
  \hline
\end{tabular}
\caption{$n=10,m=4$ (1)}\label{t4}
\end{center}

\end{table}

\begin{table}
\begin{center}
\begin{tabular}{|c|c|c|c|c|c|c|c|c|c|c|c|c|c|c|c|c|c|c|c|c|c|c|c|c|c|c|c|c|c|c|c|c|}
  \hline
  v&1&1&1&1&1&1&1&1&1&1&1 \\
  \hline
    p&1  3  8  9&1  2  7&2  3  6&2  6  7&7   9  10&1  3&1  5&1  8&1  9&3  8&3  9 \\
  \hline
  v&1&1&1.3161&1.3161&1.3161&1.1892&1.3161&1.1892&1.3161&1&\\
  \hline
  p&4  9&8  9&1&2&3&6&7&8&9&10&\\
  \hline
\end{tabular}
\caption{$n=10,m=4$ (2)}\label{t5}
\end{center}

\end{table}

\begin{table}
\begin{center}
\begin{tabular}{|c|c|c|c|c|c|c|c|c|c|c|c|c|c|c|c|c|c|c|c|c|c|c|c|c|c|c|c|c|c|c|c|c|}
  \hline
  v&1&1&1&1&1&1&1&1.3161&1.1892&1 \\
  \hline
    p&1  5  6  8&1  5  6  9&1   5   9  10&2  5  6  9&3  4  9&3   9  10&5  8  9&1  5&1  6&1  8\\
  \hline
  v&1.1892&1&1&1&1&1.3161&1&1.3161&1&1\\
  \hline
  p&1  9&1  10&2  5&2  6&2  9&5  6&5  8&5  9&5  10&6  8\\
  \hline
  v&1.1892&1&1.5651&1.1892&1.1892&1&1.7321&1.5651&1.3161&1.7321\\
  \hline
  p&6  9&9  10&1&2&3&4&5&6&8&9\\
  \hline
  v&1.3161&&&&&&&&&\\
  \hline
  p&10&&&&&&&&&\\
  \hline
\end{tabular}
\caption{$n=10,m=4$ (3)}\label{t6}
\end{center}

\end{table}
\end{example}

\section{Further Remarks}

In this paper, we studied various properties of CP tensors, showed
that a strongly symmetric, hierarchically dominated nonnegative
tensor is a CP tensor, and presented a hierarchical elimination
algorithm for checking this.    These indicate that a rich theory
for CP tensors can be established parallel to the theory of CP
matrices \cite{BS, HN, Ka, Xu}.   This theory will be a solid
foundation for applications of nonnegative tensor factorization
\cite{CZPA, SH}.   Further research on topics such as CP ranks are
needed.


\end{document}